\documentclass[a4paper,11pt]{article}
\usepackage[T1]{fontenc}
\usepackage{lmodern}
\usepackage[margin=23mm,top=20mm,bottom=20mm]{geometry}
\usepackage{amsmath,amssymb,amsthm,mathtools}
\usepackage{microtype}
\usepackage{needspace}
\usepackage{tikz}
\usepackage{float}
\usepackage{enumitem}
\usepackage[hidelinks]{hyperref}
\hypersetup{pdftitle={Triviality of metric fibrations over locally finite metric graphs},
  pdfauthor={Tsubasa Kamiyama}}
\numberwithin{equation}{section}
\newtheorem{theorem}{Theorem}[section]
\newtheorem{lemma}[theorem]{Lemma}

\theoremstyle{definition}
\newtheorem{definition}[theorem]{Definition}
\newtheorem{remark}[theorem]{Remark}
\newtheorem{example}[theorem]{Example}
\newtheorem{question}[theorem]{Question}
\newcommand{\id}{\mathrm{id}}
\newcommand{\R}{\mathbb{R}}
\newcommand{\pim}{\pi_1^{\mathrm m}}
\newcommand{\gam}[2]{\gamma^{#1}_{#2}}

\title{Triviality of metric fibrations\\over locally finite metric graphs}
\author{Tsubasa Kamiyama}
\date{November 10, 2026}

\begin{document}
\maketitle

\begin{abstract}
In this note, we give a direct proof that every metric fibration over a geodesic
metric space locally modeled on finite metric stars is isometrically isomorphic,
over its base, to the $\ell^1$ product of its base and a fiber.
Neither compactness nor completeness of the base is assumed.
We also construct a nontrivial example over a closed interval with a truncated
metric, showing that the local star condition alone does not imply the conclusion.
\end{abstract}

\section{Introduction}

The main result of this note was previously announced by the author at
Magnitude 2023 \cite{Leinster2023}.

Metric fibrations were introduced by Leinster \cite{Leinster2013} as a structure
giving a product formula for the magnitude of metric spaces.
A basic example is the triple $(X,X\times_1 F,\operatorname{pr}_X)$, where
$X\times_1 F$ is the $\ell^1$ product of the base $X$ and a metric space $F$
(Definition~\ref{def:l1-product}).
We consider conditions on the base under which every metric fibration is trivial.

Asao \cite{Asao2025} developed a classification theory of metric fibrations.
Here we construct a trivialization directly from the transport maps.

We explicitly assume that the base is geodesic and give a proof that also covers
the noncompact case.

\paragraph{Higher-dimensional bases.}
The one-dimensional case considered here leads to the following question.

\begin{question}\label{question:higher-dimensional}
For each integer $n\geq2$, what geometric conditions on an $n$-dimensional
geodesic metric space $X$ ensure that every metric fibration $(X,E,p)$ is trivial?
\end{question}

\section{Preliminaries}

\begin{definition}[Metric space]\label{def:metric}
A \emph{metric space} is a pair $(X,d_X)$ consisting of a set $X$ and a map
$d_X:X\times X\to\R_{\geq0}$ such that, for all $x,y,z\in X$,
\begin{enumerate}[label=(\arabic*),nosep,leftmargin=2.5em]
\item $d_X(x,y)=0\ \Longleftrightarrow\ x=y$;
\item $d_X(x,y)=d_X(y,x)$;
\item $d_X(x,z)\leq d_X(x,y)+d_X(y,z)$.
\end{enumerate}
All distances in this note are finite.
We write $X$ for $(X,d_X)$ when the metric is clear from the context.
\end{definition}

\Needspace{7\baselineskip}
\begin{definition}[Isometric map]\label{def:isometry}
A map $f:X\to Y$ between metric spaces $(X,d_X)$ and $(Y,d_Y)$ is
\emph{isometric} if
\[
d_Y(f(x),f(y))=d_X(x,y)\qquad(x,y\in X).
\]
\end{definition}

\begin{definition}[Isometric isomorphism]\label{def:isometric-isomorphism}
A surjective isometric map is called an \emph{isometric isomorphism}.
\end{definition}

\Needspace{7\baselineskip}
\begin{definition}[Geodesic metric space]\label{def:geodesic}
A metric space $(X,d_X)$ is \emph{geodesic} if, for every $x,y\in X$, there is
a map $c:[0,d_X(x,y)]\to X$ satisfying
\[
c(0)=x,\qquad c(d_X(x,y))=y,\qquad
d_X(c(s),c(t))=|s-t|
\]
for all $s,t\in[0,d_X(x,y)]$.
Its image is called a \emph{geodesic segment} joining $x$ to $y$.
\end{definition}

\begin{definition}[Metric star model]\label{def:star-model}
For an integer $n\geq1$ and $r>0$, set
\[
v_k=\left(\cos\frac{2\pi k}{n},\,\sin\frac{2\pi k}{n}\right)
\quad(k=0,\ldots,n-1),\qquad
S_n(r)=\bigcup_{k=0}^{n-1}\{t v_k\mid 0\leq t<r\}\subset\R^2
\]
and equip $S_n(r)$ with the metric
\[
d_{S_n(r)}(s v_i,t v_j)=
\begin{cases}
|s-t| & (i=j),\\
s+t & (i\ne j)
\end{cases}
\qquad(0\leq s,t<r).
\]
We also set $S_0(r)=\{0\}\subset\R^2$, with the zero metric.
In each case, the origin is called the \emph{center}.
\end{definition}

\begin{definition}[$n$-star neighborhood]\label{def:star-neighborhood}
Let $X$ be a metric space, $x\in X$, and $n\in\mathbb{Z}_{\geq0}$.
A subset $U\subset X$ is an \emph{$n$-star neighborhood} of $x$ if there is
a pair $(r,\varphi)$ satisfying the following conditions:
\begin{enumerate}[label=(\arabic*),nosep,leftmargin=2.5em]
\item $r\in\R_{>0}$;
\item $U$ is an open neighborhood of $x$ in $X$;
\item $\varphi:S_n(r)\to U$ is an isometric isomorphism;
\item $\varphi(0)=x$.
\end{enumerate}
Here $U$ carries the metric induced from $X$.
\end{definition}

\begin{figure}[H]
\centering
\begin{tikzpicture}[x=1cm,y=1cm,line cap=round,line join=round,
  endpoint/.style={circle,draw,fill=white,inner sep=1.6pt},
  centerpoint/.style={circle,fill,inner sep=1.5pt}]
  \node at (0,1.05) {$S_3(r)$};
  \foreach \angle in {0,120,240} {
    \draw[thick] (0,0) -- (\angle:0.8);
    \node[endpoint] at (\angle:0.8) {};
  }
  \node[centerpoint,label=below right:$0$] at (0,0) {};
  \node at (0.42,0.17) {$r$};
  \draw[->,thick] (1.25,0) -- node[above] {$\varphi$} (2.55,0);
  \draw[gray,semithick] (4.5,0.15)
    .. controls (5.55,0.15) and (5.65,1.25) .. (4.8,1.5)
    coordinate[pos=0.5] (upperright)
    .. controls (3.85,1.85) and (2.7,1.6) .. (3.1,0.8)
    coordinate[pos=0.45] (upperleft);
  \draw[gray,semithick] (4.5,0.15)
    .. controls (5.55,-0.1) and (5.6,-1.5) .. (4.55,-1.55)
    coordinate[pos=0.5] (lowerright)
    .. controls (3.5,-1.6) and (2.7,-1.4) .. (3.05,-0.8)
    coordinate[pos=0.5] (lowerleft);
  \draw[gray,semithick] (upperright)
    .. controls (5.7,0.65) and (5.9,0.75) .. (6.05,1.0);
  \draw[gray,semithick] (upperleft)
    .. controls (3.55,1.85) and (3.55,2.0) .. (3.3,2.1);
  \draw[gray,semithick] (lowerright)
    .. controls (5.55,-0.65) and (5.8,-0.75) .. (5.95,-1.0);
  \draw[gray,semithick] (lowerleft)
    .. controls (3.05,-1.6) and (3.15,-1.8) .. (2.95,-1.95);
  \draw[thick] (3.5,0) .. controls (3.85,-0.1) and (4.1,0.2) .. (4.5,0.15);
  \draw[thick] (3.5,0) .. controls (3.25,0.3) and (3.3,0.5) .. (3.1,0.8);
  \draw[thick] (3.5,0) .. controls (3.1,-0.15) and (3.2,-0.5) .. (3.05,-0.8);
  \foreach \point in {(4.5,0.15),(3.1,0.8),(3.05,-0.8)}
    \node[endpoint] at \point {};
  \node[centerpoint,label=below right:$x$] at (3.5,0) {};
  \node at (4.05,0.65) {$U$};
  \node[gray] at (5.6,1.25) {$X$};
\end{tikzpicture}
\caption{A $3$-star neighborhood $U\subset X$ with $\varphi(0)=x$ (schematic).
Open endpoints are excluded from $S_3(r)$ and $U$, respectively.}
\label{fig:star-neighborhood}
\end{figure}
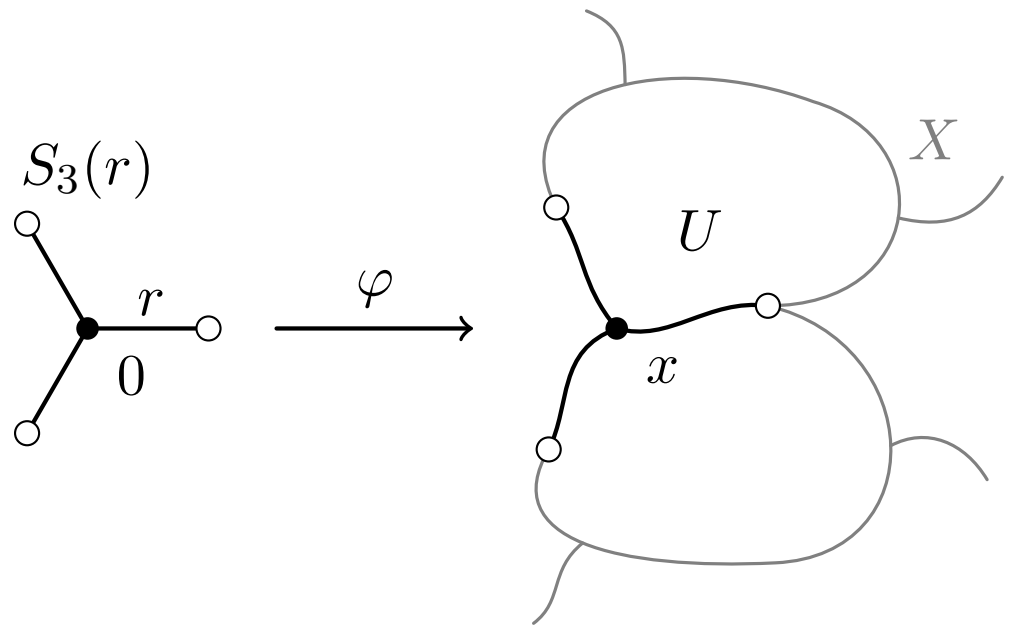

\begin{definition}[Local star structure and degree]\label{def:star}
A metric space $X$ is \emph{locally modeled on finite metric stars} if every
$x\in X$ has an $n$-star neighborhood for some integer $n\geq0$.
We call $n$ the \emph{degree} of $x$ and write $\deg_X(x)=n$.
This integer is uniquely determined: for all sufficiently small $\varepsilon>0$,
it is the number of connected components of
$\{y\in X\mid 0<d_X(x,y)<\varepsilon\}$.
\end{definition}

\begin{definition}[$1$-Lipschitz map]\label{def:lipschitz}
A map $f:X\to Y$ between metric spaces $(X,d_X)$ and $(Y,d_Y)$ is
\emph{$1$-Lipschitz} if
\[
d_Y(f(x),f(y))\leq d_X(x,y)\qquad(x,y\in X).
\]
\end{definition}

\begin{definition}[$\ell^1$ product]\label{def:l1-product}
The \emph{$\ell^1$ product} of metric spaces $(X,d_X)$ and $(F,d_F)$ is the
Cartesian product $X\times F$ equipped with the metric
\[
d_{X\times_1 F}\bigl((x,a),(y,b)\bigr)=d_X(x,y)+d_F(a,b)
\qquad(x,y\in X,\ a,b\in F).
\]
We denote this metric space by $X\times_1 F$ and call
$\operatorname{pr}_X:X\times_1 F\to X;\ (x,a)\mapsto x$ the projection.
\end{definition}

\Needspace{10\baselineskip}
\begin{definition}[Metric fibration; \cite{Leinster2013}]\label{def:fibration}
A triple $(X,E,p)$ consisting of metric spaces $(X,d_X)$ and $(E,d_E)$ and
a $1$-Lipschitz map $p:E\to X$ is a \emph{metric fibration} if, for every
$a\in E$ and $y\in X$, there exists $a_y\in p^{-1}(y)$ such that
\begin{equation}\label{eq:gate}
d_E(a,b)=d_X(p(a),y)+d_E(a_y,b)
\qquad (b\in p^{-1}(y)).
\end{equation}
We call $X$ the \emph{base}, $E$ the \emph{total space}, $p$ the
\emph{projection}, and $E_x=p^{-1}(x)$ the \emph{fiber} over $x\in X$.
Each fiber carries the metric induced from $E$.
\end{definition}

\begin{remark}\label{rem:unique-lift}
The point $a_y$ in Definition~\ref{def:fibration} is unique.
Indeed, if $a_y$ and $a'_y$ both satisfy the condition, substituting them for $b$
in \eqref{eq:gate} gives $d_E(a_y,a'_y)=0$.
\end{remark}

A metric fibration $(X,E,p)$ is \emph{trivial} if there are a metric space $F$
and an isometric isomorphism $\Psi:X\times_1 F\to E$ satisfying
$p\circ\Psi=\operatorname{pr}_X$.
Such a map $\Psi$ is called an isometric isomorphism over $X$.

\section{Canonical transport and a criterion for triviality}

Let $(X,E,p)$ be a metric fibration.

\begin{definition}[Canonical transport]\label{def:transport}
For $x,y\in X$, the map
\[
\gam{x}{y}:E_x\to E_y;\quad a\mapsto a_y
\]
is called the \emph{canonical transport} from $x$ to $y$.
Here $a_y$ is the point in Definition~\ref{def:fibration}.
\end{definition}

\begin{lemma}[{\cite[Lemma 3.4, Proposition 3.8]{Asao2025}}]\label{lem:transport}
Canonical transports are bijective isometric maps satisfying
$\gam{x}{x}=\id$ and $\gam{y}{x}=(\gam{x}{y})^{-1}$.
Moreover, for all $x,y,z\in X$ and $a\in E_x$,
\begin{equation}\label{eq:defect}
d_{E_z}\bigl(\gam{y}{z}\gam{x}{y}(a),\gam{x}{z}(a)\bigr)
\leq d_X(x,y)+d_X(y,z)-d_X(x,z).
\end{equation}
In particular, if $d_X(x,z)=d_X(x,y)+d_X(y,z)$, then
\begin{equation}\label{eq:geodesic}
\gam{y}{z}\gam{x}{y}=\gam{x}{z}.
\end{equation}
\end{lemma}

\begin{proof}
Equation~\eqref{eq:gate} gives $d_E(a,\gam{x}{y}(a))=d_X(x,y)$.
The same equation and symmetry imply
$\gam{y}{x}\gam{x}{y}(a)=a$.
For $a,b\in E_x$,
\[
d_X(x,y)+d_{E_y}(\gam{x}{y}(a),\gam{x}{y}(b))
=d_X(x,y)+d_{E_x}(a,b),
\]
so $\gam{x}{y}$ is isometric.
Put $c=\gam{y}{z}\gam{x}{y}(a)$. Then
\[
d_{E_z}(\gam{x}{z}(a),c)
=d_E(a,c)-d_X(x,z)
\leq d_X(x,y)+d_X(y,z)-d_X(x,z).
\]
Equations~\eqref{eq:defect} and \eqref{eq:geodesic} follow.
\end{proof}

Fix a basepoint $o\in X$, and set $F=E_o$ and
$\Phi_x=\gam{o}{x}:F\to E_x$.

\begin{lemma}\label{lem:criterion}
If
\begin{equation}\label{eq:criterion}
\gam{x}{y}\Phi_x=\Phi_y
\end{equation}
for all $x,y\in X$, then $\Psi(x,a)=\Phi_x(a)$ defines an isometric isomorphism
$X\times_1 F\to E$ over $X$.
Conversely, if $(X,E,p)$ is trivial, then its canonical transports satisfy
$\gam{y}{z}\gam{x}{y}=\gam{x}{z}$ for every triple of points.
\end{lemma}

\begin{proof}
The map $\Psi$ is bijective on each fiber, and \eqref{eq:criterion} gives
\[
d_E(\Phi_x(a),\Phi_y(b))
=d_X(x,y)+d_{E_y}(\Phi_y(a),\Phi_y(b))
=d_X(x,y)+d_F(a,b).
\]
For the converse, the canonical transport in a product is
$(x,a)\mapsto(y,a)$, and isometric isomorphisms over the base preserve
\eqref{eq:gate}.
\end{proof}

\section{The trivialization theorem}

\begin{lemma}\label{lem:degree-set}
Let $X$ be a metric space locally modeled on finite metric stars.
Then $B=\{x\in X\mid\deg_X(x)\ne2\}$ is closed and discrete.
In particular, $B\cap C$ is finite for every compact subset $C\subset X$.
\end{lemma}

\begin{proof}
For each $x\in X$, choose an $n$-star neighborhood $U$ of $x$.
Every point of $U\setminus\{x\}$ has a $2$-star neighborhood, so
$U\cap B\subset\{x\}$.
Thus $B$ is closed and discrete.
If $C$ is compact, then $B\cap C$ is a compact discrete space and is therefore finite.
\end{proof}

\begin{lemma}\label{lem:network}\upshape
Let $X$ be a metric space locally modeled on finite metric stars, and let $K$
be a connected union of finitely many geodesic segments containing at least two
points. Then $K$ admits a subdivision as a finite graph such that each closed
edge is contained in one of the given geodesic segments.
Each open edge is open in $X$, and each closed edge is isometric to an interval
$[0,L]$.
All endpoints of the given geodesic segments can be included among the vertices.
\end{lemma}

\begin{proof}
Put $B=\{x\in X\mid\deg_X(x)\ne2\}$.
Since $K$ is compact, Lemma~\ref{lem:degree-set} implies that $B\cap K$ is finite.
Let $V$ be the union of $B\cap K$ and all endpoints of the given geodesic segments.

Subdivide each segment at the points of $V$, obtaining finitely many isometric
images $J$ of closed intervals of positive length.
Let $J^\circ$ denote the image of the corresponding open interval.
Every point of $J^\circ$ has degree $2$ and has a sufficiently small interval
neighborhood contained in $J^\circ$.
Thus $J^\circ$ is open in $X$.
Since $J$ is compact and $J\setminus J^\circ\subset V$, the set $J^\circ$
is also closed in $K\setminus V$.
Being connected, $J^\circ$ is therefore a connected component of $K\setminus V$.
The distinct components give finitely many open edges.
Each closed edge lies in one of the given geodesic segments and has two
distinct endpoints.
\end{proof}

\begin{lemma}\label{lem:edge}
Let $\iota:[0,L]\to X;\ t\mapsto x_t$ be an isometric map, and put
$u=x_0$, $v=x_L$.
Suppose that $\iota((0,L))$ is open in $X$ and $o\notin\iota((0,L))$.
If $X$ is geodesic, then, for all $s,t\in[0,L]$,
\begin{equation}\label{eq:edge}
\gam{x_s}{x_t}=\Phi_{x_t}\Phi_{x_s}^{-1}.
\end{equation}
\end{lemma}

\begin{proof}
Put $A=d_X(o,u)$ and $B=d_X(o,v)$.
Any geodesic from $o$ to an interior point of the edge passes through $u$ or $v$,
since the open edge is open in $X$ and its boundary is contained in $\{u,v\}$.
Hence
\begin{equation}\label{eq:min}
d_X(o,x_t)=\min\{A+t,\ B+L-t\}\qquad(0\leq t\leq L).
\end{equation}
The triangle inequality gives $|A-B|\leq L$. Thus
\[
t_*=(B+L-A)/2\in[0,L].
\]
For $w=x_{t_*}$, equation~\eqref{eq:min} gives
\begin{equation}\label{eq:meeting}
\gam{u}{w}\Phi_u=\Phi_w=\gam{v}{w}\Phi_v.
\end{equation}
Equation~\eqref{eq:geodesic} and the inverse relations imply
$\gam{u}{w}=\gam{v}{w}\gam{u}{v}$.
Composing \eqref{eq:meeting} on the left with $\gam{w}{v}$ gives
\begin{equation}\label{eq:endpoints}
\Phi_v=\gam{u}{v}\Phi_u.
\end{equation}
If $t\leq t_*$, then $d_X(o,x_t)=A+t$, so
$\Phi_{x_t}=\gam{u}{x_t}\Phi_u$.
If $t\geq t_*$, the equality through $v$ and \eqref{eq:endpoints} give
\[
\Phi_{x_t}=\gam{v}{x_t}\Phi_v
=\gam{v}{x_t}\gam{u}{v}\Phi_u
=\gam{u}{x_t}\Phi_u.
\]
Equation~\eqref{eq:geodesic} and the inverse relations give
$\gam{x_s}{x_t}\gam{u}{x_s}=\gam{u}{x_t}$, which proves \eqref{eq:edge}.
\end{proof}

\begin{theorem}\label{thm:main}
Let $X$ be a nonempty geodesic metric space locally modeled on finite metric stars.
For every metric fibration $(X,E,p)$ and every basepoint $o\in X$, setting
$F=p^{-1}(o)$ gives an isometric isomorphism $X\times_1 F\cong E$ over $X$.
\end{theorem}

\begin{proof}
If $E=\varnothing$, take $F=\varnothing$.
Assume henceforth that $E\ne\varnothing$.
The definition then implies that every fiber is nonempty.
If $X=\{o\}$, the map $\Psi(o,a)=a$ gives a trivialization.
Assume that $X$ contains at least two points.
For distinct $x,y\in X$, choose geodesic segments joining $o$ to $x$, $o$ to $y$,
and $x$ to $y$, and let $K$ be their union.
A segment joining a point to itself is taken to be a singleton.

Apply Lemma~\ref{lem:network} to subdivide $K$ with $o,x,y$ among its vertices.
Lemma~\ref{lem:edge} applies to each edge.
Write $x=x_0,x_1,\ldots,x_m=y$ for the subdivision of the chosen geodesic
segment from $x$ to $y$. Then
\[
\begin{aligned}
\gam{x}{y}
&=\gam{x_{m-1}}{x_m}\cdots\gam{x_0}{x_1}\\
&=(\Phi_{x_m}\Phi_{x_{m-1}}^{-1})\cdots
  (\Phi_{x_1}\Phi_{x_0}^{-1})
=\Phi_y\Phi_x^{-1}.
\end{aligned}
\]
The first equality follows from the geodesic composition law \eqref{eq:geodesic}.
Thus \eqref{eq:criterion} holds, and Lemma~\ref{lem:criterion} gives the
trivialization $\Psi(x,a)=\gam{o}{x}(a)$.
For $x=y$, equation~\eqref{eq:criterion} follows from $\gam{x}{x}=\id$.
\end{proof}

In particular, the theorem applies to geometric realizations of finite metric
graphs equipped with the shortest-path metric.
These spaces should be distinguished from the vertex sets of graphs equipped
with the shortest-path metric.
For the latter, nontrivial examples exist over odd cycles
\cite[Proposition 5.16]{Asao2025}.

\begin{remark}\label{rem:metric-fundamental-group}
By Asao's classification \cite{Asao2025}, Theorem~\ref{thm:main} implies that
the metric fundamental group is trivial: $\pim(X,o)=\{e\}$.
\end{remark}

\Needspace{12\baselineskip}
\section{Local star structure alone does not suffice}

\begin{example}\label{ex:counterexample}
Equip $X=[0,2]$ with the metric $d_X(x,y)=\min\{|x-y|,1\}$.
This metric induces the usual topology, and each point has an open neighborhood
isometric to an interval or a half-interval.
Thus $X$ is connected, compact, and locally modeled on finite metric stars.
However, $X$ is not geodesic.

Set $h(t)=\max\{-1,\min\{t,1\}\}$ and equip $E=X\times\R$ with
\begin{equation}\label{eq:countermetric}
d_E((x,a),(y,b))=d_X(x,y)+|b-a-h(y-x)|.
\end{equation}
The function $d_E$ is symmetric and positive definite.
For every $x,y,z\in X$,
\begin{equation}\label{eq:clamp}
|h(y-x)+h(z-y)-h(z-x)|
\leq d_X(x,y)+d_X(y,z)-d_X(x,z).
\end{equation}
Indeed, for $r,s,t\in X$ with $r\leq s\leq t$, the left-hand side equals
$d_X(r,s)+d_X(s,t)-d_X(r,t)$.
Since $d_X(r,t)\geq\max\{d_X(r,s),d_X(s,t)\}$,
\[
0\leq d_X(r,s)+d_X(s,t)-d_X(r,t)
=\min\left\{
\begin{aligned}
&d_X(r,s)+d_X(s,t)-d_X(r,t),\\
&d_X(s,t)+d_X(r,t)-d_X(r,s),\\
&d_X(r,t)+d_X(r,s)-d_X(s,t)
\end{aligned}
\right\}.
\]
The left-hand side of \eqref{eq:clamp} is invariant under permutations of the
three points, so \eqref{eq:clamp} follows.
Applying the triangle inequality for absolute values and \eqref{eq:clamp}
to \eqref{eq:countermetric} gives the triangle inequality for $d_E$.

\Needspace{5\baselineskip}
For the projection $p(x,a)=x$, the triple $(X,E,p)$ is a metric fibration,
with canonical transport
\[
\gam{x}{y}:E_x\to E_y;\quad(x,a)\mapsto(y,a+h(y-x)).
\]
Equation~\eqref{eq:gate} follows from \eqref{eq:countermetric}.
Identifying each fiber with $\R$, we have
\[
\gam{1}{2}\gam{0}{1}(a)=a+2
\ne a+1=\gam{0}{2}(a).
\]
By Lemma~\ref{lem:criterion}, the fibration $(X,E,p)$ is nontrivial.
\end{example}

Thus the local star condition alone does not imply the conclusion of
Theorem~\ref{thm:main}.

\section*{Acknowledgments}

Generative AI was used to draft and revise the Japanese manuscript, to translate
it into English, and to prepare the LaTeX source.
The author is responsible for the content of this paper.

I thank Sho Shimoyama for helping me prepare my talk over burgers from MOS Burger
the day before the talk. He seems to have moved far ahead of me.
I have not had MOS Burger since that day.

\end{document}